\documentclass[12pt]{article}
\usepackage{url}
\usepackage{epic, eepic, subfigure, floatflt}
\usepackage{amsfonts, epsfig, latexsym, amsmath}
\def\rmdj {d\llap{\raise 1.22ex\hbox
  {\vrule height 0.09ex width 0.315em}\kern 0.04em}}
\def\sldj {d\llap{\raise 1.22ex\hbox
  {\vrule height 0.09ex width 0.265em}}\rlap{\raise 1.22ex\hbox
  {\vrule height 0.09ex width 0.05em}}}
\def\itdj {d\llap{\raise 1.22ex\hbox
  {\vrule height 0.09ex width 0.2em}}\rlap{\raise 1.22ex\hbox
  {\vrule height 0.09ex width 0.06em}}}
\def\bfdj {d\llap{\raise 1.16ex\hbox
  {\vrule height 0.126ex width 0.308em}\kern 0.04em}}
\def\ttdj {\rlap{\kern 0.17em\raise 1.1ex\hbox
  {\vrule height 0.09ex width 0.295em}}d}
\def\scdj {\rlap{\kern 0.04em\raise 0.57ex\hbox
  {\vrule height 0.09ex width 0.20em}}d}
\def\sfdj {d\llap{\raise 1.22ex\hbox
  {\vrule height 0.10ex width 0.3em}\kern 0.02em}}

\def\dj{\ifcase\fam \rmdj \or \or \or
  \or \itdj \or \sldj \or \bfdj \or \ttdj \or \sfdj \or \scdj \else \rmdj \fi}

\def\rmDj {\rlap{\kern 0.05em\raise 0.76ex\hbox
  {\vrule height 0.10ex width 0.28em}}D}
\def\slDj {\rlap{\kern 0.1em\raise 0.76ex\hbox
  {\vrule height 0.1ex width 0.28em}}D}
\def\itDj {\rlap{\kern 0.145em\raise 0.76ex\hbox
  {\vrule height 0.1ex width 0.274em}}D}
\def\bfDj {\rlap{\kern 0.044em\raise 0.72ex\hbox
  {\vrule height 0.126ex width 0.287em}}D}
\def\ttDj {\rlap{\kern 0.02em\raise 0.67ex\hbox
  {\vrule height 0.105ex width 0.20em}}D}
\def\scDj {\rlap{\kern 0.08em\raise 0.73ex\hbox
  {\vrule height 0.12ex width 0.24em}}D}
\def\sfDj {\rlap{\kern 0.02em\raise 0.727ex\hbox
  {\vrule height 0.126ex width 0.26em}}D}

\def\Dj{\ifcase\fam \rmDj \or \or \or
  \or \itDj \or \slDj \or \bfDj \or \ttDj \or \sfDj \or \scDj \else \rmDj \fi}

\def\QuotS#1#2{\leavevmode\kern-.0em\raise.2ex\hbox{$#1$}\kern-.1em/\kern-.1em\lower.25ex\hbox{$#2$}}

\title{\Large{\bf Classification of eight-dimensional perfect forms}}

\author{Mathieu DUTOUR SIKIRI\'C\\
Institut Ru\dj er Bo\u skovi\'c, Zagreb, Croatia\\
\ Achill SCH\"URMANN$^*$\\
Otto-von-Guericke-University, Magdeburg, Germany\\
\ Frank VALLENTIN\thanks{The second and the third author were supported by the
Deutsche Forschungsgemeinschaft (DFG) under grant SCHU 1503/4-1.
During the work on this paper the third author was also partially supported
by the Edmund Landau Center for Research in Mathematical Analysis and
Related Areas, sponsored by the Minerva Foundation (Germany), and he was
partially supported by the Netherlands Organization for Scientific
Research under grant NWO 639.032.203.}\\
CWI, Amsterdam, The Netherlands
}

\DeclareMathOperator{\Min}{Min}
\DeclareMathOperator{\Dom}{Dom}
\DeclareMathOperator{\vol}{vol}
\DeclareMathOperator{\Aut}{Aut}
\DeclareMathOperator{\Stab}{Stab}
\DeclareMathOperator{\Sym}{Sym}
\DeclareMathOperator{\Grp}{Grp}
\DeclareMathOperator{\GL}{GL}
\DeclareMathOperator{\Id}{Id}
\newcommand{\MC}{{\cal C}}
\newcommand{\MR}{{\cal M}}
\newcommand{\MF}{{\cal F}}

\newcommand{\RR}{\ensuremath{\mathbb{R}}}

\newcommand{\ZZ}{\ensuremath{\mathbb{Z}}}

\newtheorem{theorem}{Theorem}[section]

\begin{document}
\maketitle

\begin{abstract}
In this paper, we classify the perfect lattices in dimension $8$.
There are $10916$ of them.
Our classification heavily relies on exploiting symmetry in polyhedral computations.
Here we describe algorithms making the classification possible.
\end{abstract}

\section{Introduction}

A \textit{lattice} $L$ is a rank $d$ subgroup of $\RR^d$. It is of the form
$L=\ZZ v_1+\dots+\ZZ v_d$ for linearly independent $v_1, \dots, v_d$ which
are called a \textit{lattice basis} of $L$.
The determinant $\det L=\vert\det(v_1,\dots, v_d)\vert$ is independent
of the chosen basis.
By $B^d$ we denote the Euclidean unit ball.

The \textit{packing radius} $\lambda(L)$ of $L$ is defined as half
of the minimum distance between distinct lattice points.
We call $L+\lambda(L) B^d$ the \textit{lattice packing defined by $L$}.
The \textit{packing density} of the lattice packing defined by $L$ is
$\lambda(L)^d\vol(B^d)/\det\,L$, the proportion of space covered by balls.
The \textit{lattice packing problem} consists in finding the lattices
having highest packing density.

A quadratic form $q(x)$ over $\RR^d$ is a function defined as $q(x)={}^{t}xQx$ with $Q$ a real $(d\times d)$ symmetric matrix.
The quadratic form $q$ is called \textit{positive definite} if the corresponding matrix $Q$ is positive definite.
We denote by ${\cal S}^d$ (resp. ${\cal S}^d_{>0}$) the set of symmetric matrices (resp. positive definite matrices). In this paper, we will identify forms with their symmetric matrices. 
Two matrices $A, B\in {\cal S}^d_{>0}$ are called \textit{arithmetically equivalent} if there exists a $P\in \GL_d(\ZZ)$ such that $B={}^{t}PAP$.
Denote by $\Aut(A)$ the \textit{arithmetic automorphism group} of $A\in {\cal S}^d$, i.e., the group of all $P\in \GL_d(\ZZ)$ such that $A={}^{t}PAP$.

If one takes a basis $V=(v_1, \dots, v_d)$ of a lattice $L$, then $L$
can be reconstructed, up to orthogonal transformations, from the
Gram matrix $G_V=(\langle v_i, v_j\rangle)\in {\cal S}^d_{>0}$
with $\langle~,~\rangle$ being the standard scalar product on $\RR^d$.
Another basis $V'$ of $L$ corresponds to another Gram matrix 
$G_{V'}$, which is arithmetically equivalent to $G_V$.

In other words, isometry classes of lattices are in one-to-one
correspondence with arithmetical equivalence classes of positive
definite quadratic forms.
For computations, it is preferable to use the language of positive definite quadratic forms, which we will do in the remaining of the paper.
For more details on the correspondence between lattice properties and quadratic forms, see \cite[Chapter~1]{martinet}.

The \textit{arithmetical minimum} of a positive definite form $A$ is $\lambda(A)=\min_{v\in \ZZ^d-\{0\}} {}^{t}vAv$.
The \textit{Hermite invariant} of $A$ is defined as $\gamma(A)=\frac{\lambda(A)}{(\det\,A)^{1/d}}$. The \textit{packing density} of the lattice corresponding to $A$ is $\sqrt{\gamma(A)^d}\frac{\vol(B^d)}{2^d}$, with $\vol(B^d)$ being the volume of the unit ball $B^d$.
Hence, solving the lattice packing problem is equivalent to maximizing $\gamma$ over ${\cal S}^d_{>0}$.
Denote by $\gamma_d=\max_{A\in {\cal S}^d_{>0}} \gamma(A)$ the \textit{Hermite constant} in dimension $d$.

The set of \textit{shortest vectors} of $A$ is 
\begin{equation*}
\Min(A)=\{v\in \ZZ^d\mbox{~~}|\mbox{~~} {}^{t}vAv=\lambda(A)\}\;.
\end{equation*}
The form $A$ is called \textit{perfect} if the equations
\begin{equation*}
{}^{t}vBv=\lambda(A)\mbox{~for~every~}v\in \Min(A)
\end{equation*}
have the unique solution $B=A$ among $B\in {\cal S}^d$.
Perfect forms were introduced in \cite{zolotarev73}, and studied by many authors, for example, in \cite{zolotarev77}, \cite{voronoi}, \cite{CSperfect}, \cite{coxeter} and \cite{ryshkov1970}. See also \cite{zong} and \cite{martinet}.

A form $A$ is called \textit{eutactic} if there exist $\lambda_v>0$ such that
\begin{equation*}
A^{-1}=\sum_{v\in \Min(A)}  \lambda_v v{}^{t} v\;.
\end{equation*}
A form $A$ is called \textit{extreme} if the Hermite invariant $\gamma$ attains a local maximum at $A$. Voronoi proved (see \cite[\textsection~17]{voronoi} and \cite[Theorem~3.4.6]{martinet}) that a form is extreme if and only if it is perfect and eutactic.

Given a perfect form $A$, its \textit{perfect domain} is the polyhedral cone
\begin{equation*}
\Dom(A)=\bigg\{\sum_{v\in \Min(A)}\lambda_v v{}^{t}v\mbox{~}|\mbox{~} \lambda_v\geq 0\bigg\}\;.
\end{equation*}
Perfect domains form a face-to-face tessellation containing the cone
${\cal S}^d_{>0}$, i.e., every $Q\in {\cal S}^d_{>0}$ belongs to at least one
perfect domain and the intersection of two perfect domains $\Dom_1$,
$\Dom_2$ is a face of both.
If two perfect forms $A, B$ are arithmetically equivalent, i.e.,
$A={}^{t}PBP$ for some $P\in \GL_d(\ZZ)$, then $\Dom(A)=P\;\Dom(B)\; {}^{t}P$.
Voronoi proved (see \cite[\textsection~7]{voronoi} and \cite[Theorem~7.4.4]{martinet}) that for 
every fixed dimension $d$, the number of arithmetical inequivalent perfect
forms is finite.
Note, since perfect domains form a tessellation of ${\cal S}^d_{>0}$, their interest is larger than only sphere packing theory (see, for example, \cite{ash}, \cite{shepherd} and \cite{soule}).

Since the group $\GL_d(\ZZ)$ acts on the tessellation with perfect domains, this provides an algorithm for classifying all perfect forms in a fixed dimension. In particular one can solve the lattice packing problem using this algorithm.

Perfect forms have been classified up to dimension $7$:
In \cite{zolotarev77} Korkine and Zolotarev obtained the
classification of perfect forms up to dimension $5$ without using
Voronoi's algorithm.
The classification of perfect forms in dimension $6$, $7$ was done
by Barnes (see \cite{barnes}) and Jaquet (see \cite{jaquet}) using
Voronoi's algorithm.

\begin{table}
\begin{center}
\begin{tabular}{|c|c|c|c|c|}
 \hline
dim. &number. of        &authors &absolute maximum         &number of\\
     &perfect lattices  &        &of $\gamma$ realized by  & extreme lattices\\
 \hline
2       &1           &\cite{lagrange} &$\mathsf{A}_{2}=\mbox{hex}$   &1\\
3       &1           &\cite{gauss} &$\mathsf{A}_3=\mbox{fcc}$    &1\\
4       &2           &\cite{zolotarev72} &$\mathsf{D}_4$    &2\\
5       &3           &\cite{zolotarev77} &$\mathsf{D}_5$    &3\\
6       &7           &\cite{barnes}         &$\mathsf{E}_6$ &6\\
7       &33          &\cite{jaquet}         &$\mathsf{E}_7$  &30\\
8       &10916       &     &$\mathsf{E}_8$          &2408\\
 \hline
\end{tabular}
\end{center}
\caption{Perfect forms up to dimension $8$. In fact, the lattices realizing the maximum of $\gamma$ are root lattices (see \cite[Chapter~6]{CS} or \cite[Chapter~4]{martinet})}
\label{PerfectDim8}
\end{table}

\begin{theorem}\label{ClassiTheor}
There are $10916$ perfect forms in dimension $8$.
\end{theorem}
We prove the above theorem by implementing Voronoi's algorithm (see \cite{larmouth}).
This enumeration problem in dimension $8$ was considered by Martinet
and his school: in fact after the work of La{\"\i}hem \cite{laihem}, Baril \cite{baril}, Napias \cite{napias} and Batut and Martinet \cite{batutmartinet} a list of $10916$ perfect forms was known and our contribution consists in proving that this list is complete.
One key step of the enumeration is to prove:
\begin{theorem}
The polyhedral cone $\Dom(Q_{\mathsf{E}_8})$ has $25075566937584$ facets in
$83092$ orbits.
\end{theorem}

A direct consequence of the enumeration is:

\begin{theorem}
(Conjecture 6.6.7 of \cite{martinet})
Every perfect $8$-dimensional lattice has a basis of minimal vectors.
\end{theorem}

Using the face structure of the perfect domains we obtain:
\begin{theorem}\label{PossibleKissingNumber}
The set of possible \textit{kissing numbers} $|\Min(A)|$
for $A\in {\cal S}^8_{>0}$ is
$2\{1\dots 58$, $60$, $63 \dots 71$, $75$, $120\}$.
\end{theorem}

Using Theorem~\ref{ClassiTheor} Riener \cite{riener} classified all extreme lattices in dimension $8$:
\begin{theorem}\label{RienerResult}
There are $2408$ extreme lattices in dimension $8$.
\end{theorem}

Note that Mordell's inequality $\gamma_{d+1}\leq \gamma_d^{d/(d-1)}$
(for a proof see \cite{mordell}, \cite[Section~2.3]{zong} and \cite[Section~2.3]{martinet}) and the knowledge of $\gamma_7$ and $\mathsf{E}_8$ yields the value of $\gamma_8$.
The value of $\gamma_8$ was first computed by Blichfeldt in \cite{blichfeldt}.
The first proof that $\mathsf{E}_8$ is the unique form realizing $\gamma_8$ is in \cite{vetchinkin}.
Another proof based on analyzing the equality case in Mordell's inequality can be found in \cite[Section~6.6]{martinet}.
A completely different, computer assisted proof is in \cite{cohn}.
Our classification gives a fourth proof for the fact that $\mathsf{E}_8$
is the unique lattice attaining $\gamma_8$.

Basic algorithms necessary for implementing Voronoi's algorithm are
explained in Section~\ref{SectionBasicAlgo}; in Section~\ref{DualDesc}
some dual description algorithms used are explained.
In Section~\ref{Implementations} we explain specific implementation
details and in Section~\ref{ComputationalResult} some consequences
of our computation.

\section{Basic algorithms}\label{SectionBasicAlgo}
\subsection{Notions from polyhedral combinatorics}

By the Farkas-Minkowski-Weyl Theorem (see e.g.\ \cite[Corollary~7.1a]{schrijver}) a \textit{convex polyhedral cone} $\MC \subseteq \RR^m$ is defined either
by a finite set of \textit{generators} $\{v_1,\ldots, v_N\} \subseteq \RR^m$
or by a finite set of linear functionals $\{f_1, \ldots, f_M\}\subseteq (\RR^m)^*$:
\[
\MC =
\bigg\{\sum_{i=1}^N \lambda_i v_i \mbox{~}|\mbox{~} \lambda_i\geq 0\bigg\} =
\bigg\{x\in \RR^m \mbox{~}|\mbox{~} f_i(x)\geq 0\bigg\}.
\]
$\MC$ is called \textit{full-dimensional} if the only
vector space containing it is $\RR^m$.
$\MC$ is called \textit{pointed} if no linear subspace of
positive dimension is contained in it.

Let $\MC$ be a full-dimensional pointed convex polyhedral cone in $\RR^m$.
Given $f \in (\RR^{m})^{*}$, the inequality $f(x)\geq 0$
is said to be \textit{valid} for $\MC$ if it holds for all $x\in \MC$.
A \textit{face} of $\MC$ is a pointed polyhedral cone
$\{ x \in \MC \mbox{~}|\mbox{~} f(x) = 0 \}$,
where $f(x) \geq 0$ is a valid inequality.

A face of dimension $1$ is called an \textit{extreme ray} of $\MC$;
a face of dimension $m-1$ is called a \textit{facet} of $\MC$.
The set of faces of $\MC$ forms a partially ordered set under inclusion.
We write $F\lhd G$ if $F\subset G$ and $\dim F=\dim G-1$.
Two extreme rays of $\MC$ are said to be \textit{adjacent}
if they generate a two-dimensional face of $\MC$.
Two facets of ${\cal C}$ are said to be \textit{adjacent}
if their intersection has dimension $m - 2$.
Any $(m-2)$-dimensional face of $\MC$ is called a \textit{ridge}
and it is the intersection of exactly two facets of $\MC$.

Every minimal set of generators $\{v_1, \ldots, v_{N'}\}$ defining
$\MC$ has the property 
\[\{\RR_+ v_1,\ldots, \RR_+ v_{N'}\}=\{e \mbox{~}|\mbox{~} \mbox{$e$ extreme ray of $\MC$}\}.\]
Every minimal set of linear functionals $\{f_1, \ldots, f_{M'}\}$
defining $\MC$ has the property that $\{F_1, \ldots, F_{M'}\}$
with $F_i=\{x\in \MC \mbox{~}|\mbox{~} f_i(x)=0\}$ is the set of facets of $\MC$.
The problem of transforming a minimal set of generators into
a minimal set of linear functionals (or vice versa) is called
the \textit{dual description problem}.

\subsection{Voronoi's algorithm for classifying perfect forms}
In this section we describe Voronoi's algorithm.
It computes a complete representative system of arithmetically
inequivalent perfect forms:
\begin{flushleft}
\smallskip
\textbf{Input:} Dimension $d$.\\
\textbf{Output:} Set~$\MR$ of all inequivalent $d$-dimensional perfect forms.\\
\smallskip
$T \leftarrow \{Q_{\mathsf{A}_n}\}$.\\
$\MR \leftarrow \emptyset$.\\
\textbf{while} there is a $Q \in T$ \textbf{do}\\
\hspace{2ex} $\MR \leftarrow \MR \cup \{Q\}$.\\
\hspace{2ex} $T \leftarrow T \setminus \{Q\}$.\\
\hspace{2ex} $\MF \leftarrow \mbox{facets of $\Dom(Q)$}$.\\
\hspace{2ex} \textbf{for} $F \in \MF$ \textbf{do}\\
\hspace{2ex} \hspace{2ex} find perfect form $Q'$ with $F = \Dom(Q) \cap \Dom(Q')$.\\
\hspace{2ex} \hspace{2ex} \textbf{if} $Q'$ is not equivalent to a form in $\MR\cup T$ \textbf{then}\\
\hspace{2ex} \hspace{2ex}\hspace{2ex} $T \leftarrow T \cup \{Q'\}$.\\
\hspace{2ex} \hspace{2ex} \textbf{end if}\\
\hspace{2ex} \textbf{end for}\\
\textbf{end while}\\
\end{flushleft}
For the quadratic form $Q_{\mathsf{A}_n}$ in the above algorithm,
we may use $Q_{\mathsf{A}_n}=(q_{i,j})_{1\leq i,j\leq n}$
with $q_{i,i}=2$, $q_{i,i-1}=q_{i-1,i}=-1$ and $q_{i,j}=0$ otherwise
(see \cite[Section~6.1]{CS} or \cite[Section~4.2]{martinet}).

By Voronoi's finiteness theorem, we know that the above program will eventually finish.

The dual description part, which is used in computing the facets of $\Dom(Q)$, is computationally the most demanding part.
The special methods used for that purpose are explained in Section 3.
The computation of the adjacent domain is explained in Section \ref{SecAdjacent}
and the test of isometry of lattices in Section \ref{SecIsom}.

\subsection{Adjacent domain and shortest vector problems}\label{SecAdjacent}

In this section we describe the subalgorithm which computes the adjacent
perfect domain.
Given a positive definite form $A$ we need to solve the shortest vector
problem, i.e., compute its arithmetical minimum $\lambda(A)$ and the set
of vectors $\Min(A)$ that realize it.
The Fincke-Pohst algorithm (cf. \cite[Algorithm~2.7.7]{cohen}),
which has many implementations ({\tt sv} \cite{vallentinSV},
in GAP, in MAGMA, in PARI, etc.), does this.
In particular, solving a shortest vector problem in dimension $8$
is a routine task and takes only a fraction of a second.

If $F$ denotes a facet of $\Dom(A)$ with $A$ a perfect form, then we define
\begin{equation*}
\Min_F(A)=\{  v\in \Min(A)\mbox{~}|\mbox{~}v{}^{t}v\in F\}\;.
\end{equation*}
We then have the following algorithm:
\begin{flushleft}
\smallskip
\textbf{Input:} Perfect form $A$ and facet $F$ of $\Dom(A)$.\\
\textbf{Output:} Perfect form $A'$ with $\Dom(A')\cap \Dom(A)=F$.\\
\smallskip
$w \leftarrow$ an element of $\Min_F(A)$.\\
$U \leftarrow$ solution in ${\cal S}^d$ of ${}^{t}vUv=0$ for $v\in \Min_F(A)$ and ${}^{t}vUv=1$ for a\\
\hspace{3ex}\hspace{2ex} $v\in \Min(A)-\Min_F(A)$.\\
$\lambda\leftarrow 1$.\\
\textbf{while} $\Min(A+\lambda U)\subset \Min(A)$ \textbf{do}\\
\hspace{2ex} $\lambda\leftarrow 2\lambda$.\\
\textbf{end while}\\
\textbf{while} there is a $v_0\in \ZZ^n-\{0\}$ such that ${}^{t}v_0 (A+\lambda U) v_0 < {}^{t}w (A+\lambda U)w$ \textbf{do}\\
\hspace{2ex} $\lambda \leftarrow$ solution of ${}^{t}v_0 (A+ \lambda U) v_0 = {}^{t}w (A+\lambda U)w$.\\
\textbf{end while}\\
$A'\leftarrow A+\lambda U$.\\
\end{flushleft}
How do we test if there is a $v_0\in \ZZ^n-\{0\}$ such that ${}^{t}v_0 A' v_0 < {}^{t}w A'w$?
If $A'$ is positive definite this is done by solving a shortest vector problem.
If $A'$ is positive semidefinite but not positive definite we can find
a vector $v_0\in \ZZ^n-\{0\}$ in the kernel of $A'$, which will
satisfy ${}^{t}v_0 A' v_0 = 0$.
If $A'$ is not definite, we can find a vector $v_0\in \ZZ^n-\{0\}$ such
that ${}^{t}v_0 A' v_0 \leq 0$ by taking rational approximations of an
eigenvector corresponding to a negative eigenvalue.

Note that any $P\in \Aut(A)$ defines an action on $\Min(A)$ by $x\mapsto Px$.
This makes the group $\Aut(A)$ act on $\Min(A)$, on the facets of
$\Dom(A)$ and on the perfect domains adjacent to $\Dom(A)$.

\subsection{Isometry tests}\label{SecIsom}

In order to implement Voronoi's algorithm, we need to be able to
decide if two perfect forms $A$ and $B$ are arithmetically equivalent
or not.
We also need for the dual description algorithm explained later a way to
compute $\Aut(A)$.

If $A$ and $B$ are perfect forms, then $B={}^tPAP$ if and only
if $P\Min(A)=\Min(B)$. One implication is immediate; the other
follows from the definition of perfect forms.
Hence a possible algorithm is to compute $\Min(A)$ and $\Min(B)$ and
search for a $P\in \GL_n(\ZZ)$ such that $P\Min(A)=\Min(B)$; the search
space is finite since $\Min(A)$ and $\Min(B)$ are finite and form
a generating system of $\RR^{n}$.

For a non-perfect form the set of minimal vectors does
not characterize it. However,
a form $A\in{\mathcal S}^d$ is uniquely characterized
by the set of all vectors $v\in\ZZ^n$ with {\em norm}
$^t v A v$ less or equal to the maximum diagonal
coefficient $\max_{i=1}^d A_{i,i}$.
The program {\tt Isometry} (see \cite{carat}, \cite{pleskensouvignier}),
and its companion {\tt Aut\_Grp} for computing the
automorphism group of a form, work by generating these
possibly very large vector sets.
Hence in order to speed up computations with these
programs, it is desirable to compute with forms with small,
maximum diagonal coefficients. This can be achieved by
storing only a {\em Minkowski reduced form} for each
equivalence class. Such can be obtained for example
with the program {\tt Mink\_red} (see \cite{carat}).

\section{Dual description methods}\label{DualDesc}
General purpose programs like {\tt cdd} \cite{cdd}, {\tt lrs} \cite{lrs},
{\tt pd} \cite{pd} and {\tt porta} \cite{porta}
allow one to compute the dual description of a polyhedral cone given by its
facets (linear inequalities) or by its extreme rays (generators).
Since the programs are implementations
of quite different methods, their efficiency may vary tremendously
on a particular cone. The perfect form $Q_{\mathsf{E}_8}$ (see \cite[Chapter~6]{CS})
has $|\Min(Q_{\mathsf{E}_8})|=240$, therefore the $36$-dimensional cone
$\Dom(Q_{\mathsf{E}_8})$ has $120$ extreme rays.
All general purpose programs take too much time in computing the
facets of this cone.

Usually one is interested only in a list of representatives of orbits
of facets and not in the full list of facets.
This leads us naturally to the \textit{Adjacency Decomposition Method},
which exploits the symmetry of a polyhedral cone.
The programs mentioned above are still used, but as a subroutine.
In the following discussion, we will assume that $\MC$ is a full
dimensional, pointed polyhedral cone in $\RR^{m}$ generated by
the extreme rays $(e_i)_{1\leq i\leq N}$ and we want to compute
its facets.

\subsection{Adjacency Decomposition Method}

\begin{flushleft}
\smallskip
\textbf{Input:} Extreme rays of a polyhedral cone $\MC$ and a group $G$ acting on $\MC$.\\
\textbf{Output:} Complete set~$\MR$ of inequivalent facets of $\MC$ under $G$.\\
\smallskip
$T \leftarrow \{F\}$ with $F$ a facet of $\MC$.\\
$\MR \leftarrow \emptyset$.\\
\textbf{while} there is a $F \in T$ \textbf{do}\\
\hspace{2ex} $\MR \leftarrow \MR \cup \{F\}$.\\
\hspace{2ex} $T \leftarrow T \setminus \{F\}$.\\
\hspace{2ex} $\MF \leftarrow \mbox{facets of $F$}$.\\
\hspace{2ex} \textbf{for} $H \in \MF$ \textbf{do}\\
\hspace{2ex} \hspace{2ex} find facet $F'$ of $\MC$ with $H = F \cap F'$.\\
\hspace{2ex} \hspace{2ex} \textbf{if} $F'$ is not equivalent under $G$ to a facet in $\MR\cup T$ \textbf{then}\\
\hspace{2ex} \hspace{2ex}\hspace{2ex} $T \leftarrow T \cup \{F'\}$.\\
\hspace{2ex} \hspace{2ex} \textbf{end if}\\
\hspace{2ex} \textbf{end for}\\
\textbf{end while}\\
\end{flushleft}
Note that Voronoi's algorithm is very similar to the Adjacency
Decomposition Method. Both fit into the framework of
\textit{graph traversal algorithms}.
In the above algorithm, an initial
facet can be found by solving a linear program.

The algorithm relies on the ability to test if two facets are
equivalent under the symmetry group.
The possible strategies and the one that we used are explained in
Section \ref{Implementations}.
The Adjacency Decomposition Method framework is a reasonably natural
algorithm for computing with symmetry. Hence it was discovered several
times, for example, in \cite{jaquet} as ``algorithm de l'explorateur'', 
in \cite{CR} as ``adjacency decomposition method'' and
in \cite{DFPS} as ``subpolytope algorithm''.

We explain here the {\em gift-wrapping step} (see \cite{Swa})
to compute an adjacent facet.
Given the list of extreme rays $(e_i)_{1\leq i\leq N}$, a facet $F\lhd \MC$
is encoded by an index set $S_F\subset \{1,\dots, N\}$ such that $F$
is generated by $(e_i)_{i\in S_F}$.
Given a ridge $H$, we need an algorithm for computing the uniquely
determined facet $F'$ of $\MC$ such that $F\cap F'=H$.
The ridge $H$ is encoded by a set $S_H\subset S_F$ such that
$(e_i)_{i\in S_H}$ generate $H$.
The defining inequalities $f\in (\RR^m)^*$ of the facet $F'$
should satisfy $f(e_i)=0$ for all $i\in S_H$.
The vector space of such functions has dimension $2$.
Let us select a basis $\{f_1, f_2\}$ of it.
If $f=\alpha_1 f_1+\alpha_2 f_2$ is the defining inequality of $F$ or $F'$,
$f(e_i)\geq 0$ for all $i$ with $1\leq i\leq N$.
This translates into a set of linear inequalities on $\alpha_1, \alpha_2$
defining a $2$-dimensional pointed polyhedral cone.
One easily finds its two generators $(\alpha^{i}_1, \alpha^{i}_2)_{1\leq i\leq 2}$.
The corresponding inequalities
$f_i(x)=\alpha^i_1 f_1(x)+\alpha^i_2 f_2(x)\geq 0$ on $\MC$ define 
the two adjacent facets $F$ and $F'$ of $\MC$.
\bigskip

The Adjacency Decomposition Method can find the dual description
of very symmetric polyhedral cones, when other methods fail.
But this algorithm uses dual description, albeit in one dimension
lower and again this computation might be impossible by the known
general purpose algorithms.
The \textit{incidence number} of a face is the number of extreme
rays contained in it;
from our experience it is a good measure of the complexity of a
face: in all polyhedral cones encountered by us so far,
the facets with the highest incidence number are the ones of
highest symmetry and they
are also the ones whose dual description is the most difficult to compute.

As a consequence, we begin the computation from the orbit with the lowest
incidence number, since they are presumably easiest to treat and 
we may not need to treat all orbits, because of the following theorem
due to Balinski:

\begin{theorem}\label{Balinski}(\cite{balinski}, see e.g. \cite[Theorem~3.14]{ziegler})
Let $\MC$ be an $m$-dimensional, pointed polyhedral cone. Let $G$ be the
undirected graph whose vertices are the facets of $\MC$ and whose
edges are the ridges of $\MC$. Two vertices $E_1, E_2$
are connected by an edge $F$ if $E_1\cap E_2=F$.
Then, the graph $G$ is $(m - 1)$-connected, i.e., removal of any
$m - 2$ vertices leaves it connected.
\end{theorem}
Using the above theorem, we know that if the total number of facets in
unfinished orbits is less than $m-1$, then they cannot be adjacent
to yet to be discovered facets and so we are done.
In practice this simple criterion, which can be considered as
an extension of Th\'eor\`eme 7 and its corollaries in \cite{jaquet}, is 
extremely useful and many difficult computations were finished by it.

\subsection{Recursive Adjacency Decomposition Method}

When there are only a few remaining facets to treat and we cannot apply
Theorem~\ref{Balinski}, our method is to use the Adjacency Decomposition
Method on the remaining untreated facets of $\MC$ recursively.
The problem is that one may be confronted with a lot of cases to
consider. In this section we present the method used to make this
manageable in some cases.

If a face $F'$ satisfies $F'\lhd F_1\lhd \MC$ for 
a facet $F_1$, then there is exactly one other facet $F_2$ such that 
$F'\lhd F_2\lhd \MC$. Hence, if one applies the Adjacency Decomposition
Method to $F_1$ and $F_2$, then one will compute the dual description
of $F'$ two times.
The number of such repetitions increases as the recursion depth
increases.

To handle this, we use a \textit{banking system}.
We store the representatives of orbits of $(k-1)$-dimensional
faces (facets) $F$ of a $k$-face $F'$ with respect to a
group $\Grp_{F'}\subseteq \Aut_{F'}$ of linear automorphisms
of $F'$, which is not necessarily the full automorphism group of $F'$.

If $F'\lhd F''$ and the Adjacency Decomposition Method is applied to
$F''$, then as a subtask one needs to find the orbits of facets of $F'$ under
the stabilizer $\Stab(\Grp_{F''}, F')$ of $F'$ under the group $\Grp_{F''}$.
So, the problem is to obtain the list of facets under the action of
$\Stab(\Grp_{F''}, F')$
from a list of orbits under the action of $\Grp_{F'}$.

When $\Stab(\Grp_{F''}, F')$ is not a subgroup of $\Grp_{F'}$
(we did not require that $\Grp_{F'} = \Aut_{F'}$), then we replace
$\Grp_{F'}$ by the group generated by $\Grp_{F'}$ and $\Stab(\Grp_{F''}, F')$.
So, we can assume that $\Stab(\Grp_{F''}, F')\subseteq \Grp_{F'}$.
Take an orbit $\Grp_{F'}F$ of $(k-1)$-dimensional faces of the $k$-dimensional face $F'$.
We find elements $g_1, \dots, g_r$ in $\Grp_{F'}$ such that
\begin{equation*}
\Grp_{F'}=\bigcup_{i=1}^r \Stab(\Grp_{F''},F') g_i\Stab(\Grp_{F'},F),
\end{equation*}
i.e., we compute a decomposition of $\Grp_{F'}$ into \textit{double cosets}
(using a computer algebra system like {\tt GAP} \cite{GAP}).
One then obtains $\Grp_{F'} F=\bigcup_{i=1}^r \Stab(\Grp_{F''}, F') g_iF$,
i.e., the orbit $\Grp_{F'} F$
splits into $r$ orbits $\Stab(\Grp_{F''}, F') F_i$ with $F_i=g_iF$.
This double coset decomposition is a classic enumeration technique,
exposed for example in \cite{brinkmann} and \cite{kerber}.

To set-up a banking system as described, we need to be able to 
test if two polyhedral cones are isomorphic and we need to compute
their automorphism groups.

\subsection{Isomorphisms and automorphisms of polyhedral cones}

The symmetry group of a pointed polyhedral full dimensional cone in $\RR^m$
generated by extreme rays $(e_i)_{1\leq i\leq N}$
is the group of matrices $A\in \GL_m(\RR)$ such that there
exists a permutation $\sigma$ of $\{1,\dots,N\}$
with $Ae_i=e_{\sigma(i)}$ for $1\leq i\leq N$.
This infinite group is called \textit{projective automorphism group} (see \cite{kaibelschwartz}).

We are not aware of an algorithm to decide whether or not polyhedral
cones are equivalent and to compute their projective automorphism groups.
However, for the stronger notion of restricted equivalence introduced
below, we can resolve those questions. 
Computing with a proper subgroup of the projective automorphism
group is not a problem for the Recursive Adjacency Decomposition Method.
However, it impacts the computing time since we may compute 
the dual description of projectively isomorphic polyhedral cones several
times.

Given a full-dimensional vector family $(v_i)_{1\leq i\leq N}$,
we define the positive definite matrix
\begin{equation*}
Q=\sum_{i=1}^{N} v_i{}^{t}v_i \in {\cal S}^m_{>0}\;.
\end{equation*}
Furthermore, denote by $R$ the unique matrix $R\in {\cal S}^m_{>0}$
satisfying $Q^{-1}=R^2$ and $(w_i)_{1\leq i\leq N}=(Rv_i)_{1\leq i\leq N}$ the
image of the family $(v_i)$ under $R$.
Define $G(v_i)$ to be the complete graph with vertices $v_i$
and edge weights $c_{ij}={}^{t}v_i Q^{-1} v_j=\langle w_i, w_j\rangle$.

A \textit{restricted isomorphism} of two vector families
$(v_i)_{1\leq i\leq N}$ and $(v'_i)_{1\leq i\leq N}$ is given by
a matrix $A$ such that there exists a permutation $\sigma$ satisfying
$Av_i=v'_{\sigma(i)}$ for $i=1,\dots, N$.

Such a restricted isomorphism satisfies $AQ{}^{t}A=Q'$
with $Q$ and $Q'$ as above. One then checks that the matrix
$O=R'AR^{-1}$ is orthogonal and satisfies $Ow_i=w'_{\sigma(i)}$.
This implies $c_{ij}=c'_{\sigma(i)\sigma(j)}$; i.e., the restricted
isomorphism of the vector families $(v_i)$ and $(v'_i)$ corresponds
to an isomorphism between the edge weighted graph $G(v_i)$ and $G(v'_i)$.

Now we need to prove that if $\sigma$ is an isomorphism between $G(v_i)$ and
$G(v'_i)$, then the equation $Av_i=v'_{\sigma(i)}$ admits only one solution.
Clearly, we can assume $\sigma=\Id$ and consider the equivalent equation
$Ow_i=w'_i$.
Since the $v_i$ generate $\RR^m$, we find a basis
$(v_{i_1}, \dots, v_{i_m})$ of $\RR^m$.
If $P$, respectively $P'$, is the $m\times m$ matrix formed by $(w_{i_k})$,
respectively $(w'_{i_k})$, then the equation $c_{ij}=c'_{ij}$ takes the
form ${}^{t}PP={}^{t}P'P'$. So, the matrix $O=P'P^{-1}$ is orthogonal
and one has for any $1\leq k\leq m$, $1\leq j\leq N$:
\begin{equation*}
\langle w'_{i_k}, Ow_j\rangle=\langle Ow_{i_k}, Ow_j\rangle=\langle w_{i_k}, w_j\rangle=\langle w'_{i_k}, w'_j\rangle.
\end{equation*}
The above equation gives $\langle w'_i, Ow_j-w'_j\rangle=0$; since the
$w'_i$ form a basis of $\RR^m$ we obtain the relation $Ow_j=w'_j$,
i.e., a restricted isomorphism between $(v_i)$ and $(v'_i)$.
In the same manner one proves that the restricted automorphism problem
for a vector family $(v_i)_{1\leq i\leq N}$ is equivalent to the
automorphism problem of the edge weighted graph $G(v_i)$.

\section{Implementation details}\label{Implementations}

The key part of our algorithm is to test if
two $(k-1)$-dimensional faces of a $k$-dimensional face $F$
are equivalent under a group $\Grp_F$ of automorphisms of $F$.
In practice, we represent $(k-1)$-dimensional faces by the set
of indices of extreme rays contained in $F$ and we use the 
command {\sf RepresentativeAction} with the action {\sf OnSets}
of the GAP computer algebra system \cite{GAP} to test for equivalence.
We compute the stabilizer of a face with the {\sf Stabilizer} command.
Building the full orbit, a strategy we could not consider, is used
in \cite{DFPS} and \cite{mining}.
For some special groups and representations, like $\Sym(n)$ acting
on $n$ elements, it is very easy to test equivalence; this
strategy is used in \cite{anzin} and \cite{robbins}.

In practice, the program {\tt nauty} \cite{nauty} computes
efficiently isomorphisms and automorphisms for vertex
weighted graphs. Its ``User's Guide'' (Version 2.4, page 25)
contains a description of a simple method to transform edge
weighted problems into vertex weighted ones.

\section{Computational results}\label{ComputationalResult}

To compute the dual description of the $10916$ perfect domains
we only had to use the Adjacency Decomposition Method in two cases:
for $\mathsf{E}_8$ and for the Barnes-Coxeter lattice $\mathsf{A}_8^2$ (cf. \cite[Section~5.1]{martinet}).
It would also be necessary for the lattice $\mathsf{D}_8$
(cf. \cite{jaquet92}), but for $\mathsf{D}_n$ the list of neighboring
domains is already known; see \cite[Theorem~15]{baranovski}.
All other domains could be treated by the standard software {\tt lrs}.

Computing the dual description of $\Dom(Q_{\mathsf{E}_8})$ (whose
symmetry group is $\QuotS{W(E_8)}{\ZZ_2}$ of size $348364800$)
took several months of computer time.

Most of the time of this computation was used for treating
three orbits of facets of $\Dom(Q_{\mathsf{E}_8})$
generated by $66$, $70$ and $75$ extreme rays respectively.
The facet with $75$ extreme rays is interesting:
its stabilizer under the action of the group $\QuotS{W(\mathsf{E}_8)}{\ZZ_2}$
has size $23040$,
but when the automorphism group is computed, the size grows to
$737280$, therefore allowing us to finish the computation.
We have no theoretical explanation for those additional symmetries.

It is proved in \cite{watson} that if a positive definite form $A$
has $\frac{1}{2}\vert\Min(A)\vert > 75$, then $A$ is arithmetically equivalent
to a multiple of $Q_{\mathsf{E}_8}$.
The knowledge of facets of the perfect domain allows us to say
a bit more: If $\frac{1}{2}|\Min(A)|=75$, then $A=\lambda_1 A_1+\lambda_2 A_2$
with both $A_i$ arithmetically equivalent to $Q_{\mathsf{E}_8}$
and $\lambda_i>0$.
If $\frac{1}{2}|\Min(A)|=70$, then $A=\lambda_1 A_1 + \lambda_2 A_2$
with $\lambda_i>0$, $A_1$ arithmetically equivalent to
$Q_{\mathsf{E}_8}$ and $A_2$ arithmetically
equivalent to $Q_{\mathsf{A}_8^2}$.

It is also worthwhile to note that the form $Q_{\mathsf{E}_8}$ is
not contiguous to only
two forms: $Q_{\mathsf{A}_8}$ and the one with number 8190 (see \url{http://www.math.uni-magdeburg.de/lattice_geometry/} for the contiguities).

The method developed here allowed us to treat the $8$-dimensional case.
Without using symmetry of perfect domains, one would be limited to
dimension $6$.
Note that the algorithm developed by Jaquet in \cite{jaquet} for dimension
$7$ is also a Recursive Adjacency Decomposition Method; he does not use a
banking system and instead of using a standard software for the double description (which did not exist at that time) he uses a specific ``Cascade algorithm''.
We computed the facets of the perfect domain $\Dom(Q_{\mathsf{E}_7})$
in less than a day.

In dimension $9$, a priori it seems currently very difficult to do a full
computation of the list of perfect forms. 
Their number becomes even larger. We found more than $500000$
(see \url{http://www.math.uni-magdeburg.de/lattice_geometry/}) and there is no end in sight. Nevertheless, it is probably only a matter of time until the perfect forms in dimension $9$ can be enumerated by Voronoi's algorithm.

\section{Erratum to algorithm of section 2.3}

The algorithm given in section 2.3 is incorrect.
It does not work when implemented and the computations of the paper were
done with a different and correct algorithm, which is explained below:
\begin{flushleft}
\smallskip
\textbf{Input:} Perfect form $A$ and facet $F$ of $\Dom(A)$.\\
\textbf{Output:} Perfect form $A'$ with $\Dom(A')\cap \Dom(A)=F$.\\
\smallskip
$(l,u)\leftarrow (0,1)$\\
\textbf{while} $A+ u F\notin {\cal S}^{d}_{>0}$ or $\lambda(A + u F)=\lambda(A)$ \textbf{do}\\
\hspace{2ex}\textbf{if} $A+ u F\notin {\cal S}^d_{>0}$ \textbf{then}\\
\hspace{2ex}\hspace{2ex}$u\leftarrow (l+u)/2$\\
\hspace{2ex}\textbf{else}\\
\hspace{2ex}\hspace{2ex}$(l,u) \leftarrow (u,2u)$\\
\hspace{2ex}\textbf{end if}\\
\textbf{end while}\\
\textbf{while} $\Min(A + l F)\subseteq \Min(Q)$ \textbf{do}\\
\hspace{2ex}$\gamma\leftarrow (l+u)/2$\\
\hspace{2ex}\textbf{if} $\lambda(A + \gamma F) \geq \lambda(A)$ \textbf{then}\\
\hspace{2ex}\hspace{2ex}$l\leftarrow \gamma$\\
\hspace{2ex}\textbf{else}\\
\hspace{2ex}\hspace{2ex}$u\leftarrow \min\{(\lambda(A) - A[v])/F[v]\quad : \quad v\in \Min(A + \gamma F), F[v] < 0\}\cup \{\gamma\}$\\
\hspace{2ex}\textbf{end if}\\
\hspace{2ex}\textbf{if} $\lambda(A + u F) = \lambda(A)$ \textbf{then}\\
\hspace{2ex}\hspace{2ex}$l\leftarrow u$\\
\hspace{2ex}\textbf{end if}\\
\textbf{end while}\\
\end{flushleft}
The first section of the algorithm determines two values $l$ and $u$ which are lower and upper bounds for the value of the lifted parameter.
The second section adjust until we find a correct value of the parameters based on shortest vector computations. The arrangement of tests is subtle but guarantees a termination in any case.

\end{document}